\documentclass[a4paper,11pt]{amsart}
\usepackage[a4paper,twoside,centering]{geometry}
\usepackage{amsfonts}

\title[Hochschild cohomology of cluster-tilted algebras]
{Hochschild cohomology of the cluster-tilted algebras of finite
representation type}
\author{Sefi Ladkani}
\address{%
Mathematical Institute of the University of Bonn \\
Endenicher Allee 60 \\
53115 Bonn, Germany}
\urladdr{http://www.math.uni-bonn.de/people/sefil}
\email{sefil@math.uni-bonn.de}

\DeclareMathOperator{\ch}{char}
\DeclareMathOperator{\Ext}{Ext}
\DeclareMathOperator{\HH}{HH}

\newcommand{\eps}{\varepsilon}
\newcommand{\gL}{\Lambda}
\newcommand{\bZ}{\mathbb{Z}}

\theoremstyle{plain}
\newtheorem{theorem}{Theorem}
\newtheorem{cor}{Corollary}

\theoremstyle{definition}
\newtheorem*{remark*}{Remark}

\begin{document}
\begin{abstract}
We compute the Hochschild cohomology groups of the cluster-tilted
algebras of finite representation type.
\end{abstract}

\maketitle

An important homological invariant of a finite-dimensional algebra
$\gL$ over a field $K$ is its Hochschild cohomology, defined as the
graded ring $\HH^*(\gL) = \Ext^*_{\gL^{op} \otimes_K \gL}(\gL,\gL)$,
see~\cite{Happel89}. Even if $\gL$ is given combinatorially as quiver
with relations, it is not easy to explicitly determine the groups
$\HH^i(\gL)$, and in many cases one needs a projective resolution of
$\gL$ as a bimodule over itself.

An interesting class of algebras consists of the cluster-tilted
algebras introduced by Buan, Marsh and Reiten~\cite{BMR07} as the
endomorphism algebras of cluster-tilting objects in a cluster category.
Cluster-tilted algebras of finite representation type were studied
in~\cite{BMR06}, see also~\cite{CCS06b}. They can be described by
quivers with relations where the quivers are obtained from orientations
of $ADE$ Dynkin diagrams by performing sequences of quiver
mutations~\cite{FominZelevinsky02}, and the defining relations consist
of zero- and commutativity-relations that can be deduced from the
quiver in an algorithmic way. From a homological viewpoint,
cluster-tilted algebras are Gorenstein~\cite{KellerReiten07}, but in
general they are of infinite global dimension.

Previously, the first Hochschild cohomology group of a schurian
cluster-tilted algebra was computed in terms of an underlying tilted
algebra, see~\cite{AssemRedondo09}. In this note we compute all the
Hochschild cohomology groups $\HH^i(\gL)$ of a cluster-tilted algebra
$\gL$ of finite representation type in terms of its quiver.

In order to formulate our results, we encode the dimensions of
$\HH^i(\gL)$ in a formal power series
\[
h_{\gL}(z) = \sum_{i=0}^{\infty} \dim_K \HH^i(\gL) \cdot z^i - 1
\]
and we define, for $n \geq 3$, the formal power series
\[
f_n(z) = \frac{z}{1-z} -
\frac{z^2 \left(1 + \eps_n(z + z^2) + z^3 \right)}{1 - z^{2n}},
\quad \text{where }
\eps_n = \begin{cases}
0 & \text{if $\ch K$ divides $n-1$,} \\
1 & \text{otherwise.}
\end{cases}
\]

The cluster-tilted algebras of Dynkin type $A$ have been described as
quivers with relations in~\cite{BuanVatne08,CCS06a}.

\begin{theorem}[Dynkin type A]
Let $\gL$ be a cluster-tilted algebra of Dynkin type $A$ and let $t$ be
the number of oriented $3$-cycles in the quiver of $\gL$. Then
$h_\gL(z) = t f_3(z)$.
\end{theorem}

The quivers in the mutation class of a Dynkin quiver of type $D$ have
been explicitly described in~\cite{Vatne10}, where they were organized
into four types. In the next theorem we use the terminology
of~\cite[\S1.6]{BHLtypeD} concerning types and parameters.

\begin{theorem}[Dynkin type D]
Let $\gL$ be a cluster-tilted algebra of Dynkin type $D$.
\begin{enumerate}
\renewcommand{\theenumi}{\Roman{enumi}}
\item
If $\gL$ is of type I with parameters $(s,t)$, then
\[
h_\gL(z) = t f_3(z) .
\]

\item
If $\gL$ is of type II with parameters $(s_1,t_1,s_2,t_2)$, then
\[
h_\gL(z) = (1+t_1+t_2) f_3(z) .
\]

\item
If $\gL$ is of type III with parameters $(s_1,t_1,s_2,t_2)$, then
\[
h_{\gL}(z) = f_4(z) + (t_1+t_2)f_3(z) .
\]

\item [(IVa)]
If $\gL$ is of type IV and its quiver is an oriented cycle of length
$n$, then
\[
h_{\gL}(z) = f_n(z) .
\]

\item [(IVb)]
If $\gL$ is of type IV with parameters $\bigl((d_1,s_1,t_1),
(d_2,s_2,t_2), \dots, (d_r, s_r, t_r)\bigr)$, then
\[
h_{\gL}(z) = f_n(z) + (t_1 + t_2 + \dots + t_r) f_3(z) ,
\]
where $n = d_1 + \dots + d_r + \left|\{ 1 \leq j \leq r \,:\, d_j=1 \}
\right|$.
\end{enumerate}
\end{theorem}

\begin{remark*}
The cluster-tilted algebra corresponding to an oriented cycle is a
truncated cycle algebra. The Hochschild cohomology of such algebras was
considered by many authors,
see~\cite{BLM00,ErdmannHolm99,Locateli99,Zhang97}.
\end{remark*}

We define the \emph{associated polynomial} of an algebra $\gL$ as
$\det(xC_\gL - C^T_{\gL}) \in \bZ[x]$, where $C_{\gL}$ denotes the
Cartan matrix of $\gL$. The cluster-tilted algebras of Dynkin type $E$
have been classified up to derived equivalence in~\cite{BHLtypeE},
where it is shown that the associated polynomial is a complete derived
invariant for these algebras.

\begin{theorem}[Dynkin type E]
Let $\gL$ be a cluster-tilted algebra of Dynkin type $E$. Then
$h_\gL(z)$ is determined by the associated polynomial of $\gL$
according to Table~\ref{tab:E}.
\end{theorem}

\begin{table}
\[
\begin{array}{lc}
\multicolumn{1}{c}{\textbf{Associated polynomial}} &
\mathbf{h_{\gL}(z)} \\
\hline \\[-11pt]
x^6-x^5+x^3-x+1 & 0 \\
2(x^6-2x^4+4x^3-2x^2+1) & f_3(z) \\
2(x^6-x^4+2x^3-x^2+1) & f_3(z) \\
3(x^6+x^3+1) & f_4(z) \\
4(x^6+x^4+x^2+1) & f_5(z) \\
4(x^6+x^5-x^4+2x^3-x^2+x+1) & 2f_3(z) \\
\hline \\[-11pt]
x^7-x^6+x^4-x^3+x-1 & 0 \\
2(x^7-2x^5+4x^4-4x^3+2x^2-1) & f_3(z) \\
2(x^7-x^5+x^4-x^3+x^2-1) & f_3(z) \\
2(x^7-x^5+2x^4-2x^3+x^2-1) & f_3(z) \\
3(x^7-1) & f_4(z) \\
4(x^7+x^5-2x^4+2x^3-x^2-1) & f_5(z) \\
4(x^7+x^5-x^4+x^3-x^2-1) & f_5(z) \\
4(x^7+x^6-2x^5+2x^4-2x^3+2x^2-x-1) & 2f_3(z) \\
4(x^7+x^6-x^5-x^4+x^3+x^2-x-1) & 2f_3(z) \\
4(x^7+x^6-x^5+x^4-x^3+x^2-x-1) & 2f_3(z) \\
5(x^7+x^5-x^4+x^3-x^2-1) & f_6(z) \\
6(x^7+x^5-x^2-1) & f_7(z) \\
6(x^7+x^6-x^4+x^3-x-1) & f_4(z)+f_3(z) \\
8(x^7+x^6+x^5-x^4+x^3-x^2-x-1) & f_5(z)+f_3(z) \\
\hline \\[-11pt]
x^8-x^7+x^5-x^4+x^3-x+1 & 0 \\
2(x^8-2x^6+4x^5-4x^4+4x^3-2x^2+1) & f_3(z) \\
2(x^8-x^6+x^5+x^3-x^2+1) & f_3(z) \\
2(x^8-x^6+2x^5-2x^4+2x^3-x^2+1) & f_3(z) \\
3(x^8+x^4+1) & f_4(z) \\
4(x^8+x^6-2x^5+4x^4-2x^3+x^2+1) & f_5(z) \\
4(x^8+x^6-x^5+2x^4-x^3+x^2+1) & f_5(z) \\
4(x^8+x^7-2x^6+2x^5+2x^3-2x^2+x+1) & 2f_3(z) \\
4(x^8+x^7-x^6+2x^4-x^2+x+1) & 2f_3(z) \\
4(x^8+x^7-x^6+x^5+x^3-x^2+x+1) & 2f_3(z) \\
5(x^8+x^6+x^4+x^2+1) & f_6(z) \\
6(x^8+x^6+x^5+x^3+x^2+1) & f_7(z) \\
6(x^8+x^7+2x^4+x+1) & f_4(z)+f_3(z) \\
8(x^8+x^7+x^6+2x^4+x^2+x+1) & f_5(z)+f_3(z) \\
8(x^8+2x^7+2x^4+2x+1) & 3f_3(z)
\end{array}
\]
\caption{The Hochschild cohomology groups as functions of the
associated polynomial for cluster-tilted algebras of types $E_6, E_7,
E_8$.} \label{tab:E}
\end{table}

We list a few consequences of these results. The first states that
cluster-tilted algebras of finite representation type are rigid.

\begin{cor}
$\HH^2(\gL) = 0$ for any cluster-tilted algebra of finite
representation type $\gL$.
\end{cor}

Another consequence is that the Hochschild cohomology groups of a
cluster-tilted algebra of finite representation type are completely
determined by its first Hochschild cohomology and the determinant of
its Cartan matrix. These determinants were computed
in~\cite{BHLtypeE,BHLtypeD,BuanVatne08}.

\begin{cor}
Let $\gL$ be a cluster-tilted algebra of finite representation type.
Then
\[
h_{\gL}(z) = f_n(z) + t f_3(z)
\]
where $t = \dim \HH^1(\gL) - 1$ and $n = 1 + 2^{-t} \det C_{\gL}$.

In particular, for two cluster-tilted algebras $\gL, \gL'$ of finite
representation type the following conditions are equivalent:
\begin{enumerate}
\renewcommand{\theenumi}{\roman{enumi}}
\item
$\HH^1(\gL) \simeq \HH^1(\gL')$ and $\det C_{\gL} = \det C_{\gL'}$;

\item
$\HH^i(\gL) \simeq \HH^i(\gL')$ for all $i \geq 0$.
\end{enumerate}
\end{cor}

\begin{remark*}
Fixing the number of simples, we see that the Hochschild cohomology is
a complete derived invariant for cluster-tilted algebras of Dynkin type
$A$. This is no longer true in Dynkin types $D$ and $E$.
\end{remark*}

Our results are based on several ingredients. The first is the explicit
knowledge of the quivers of the cluster-tilted algebras in
question~\cite{BHLtypeE,BuanVatne08,CCS06a,Vatne10}. The second
ingredient is a reduction technique based on the long exact sequences
of~\cite{Cibils00,Happel89,MichelenaPlatzeck00} allowing one to
decompose the problem of computing the Hochschild cohomology of a
cluster-tilted algebra of finite representation type into smaller
problems involving simpler cluster-tilted algebras. However, some of
these simpler algebras are not monomial, hence the projective
resolution given in~\cite{Bardzell97} is not always applicable. In
order to overcome this difficulty, we use the invariance of Hochschild
cohomology under derived equivalence~\cite{Happel89,Rickard91} and
replace these algebras by derived equivalent ones whose quivers are
oriented cycles and their defining relations consist of only
zero-relations of varying lengths. In general, these monomial algebras
are not cluster-tilted anymore. Finally, by applying another reduction
technique we are able to shorten the cycles and show that the
Hochschild cohomology of these monomial algebras is isomorphic to that
of certain truncated cycle algebras. The Hochschild cohomology of the
latter algebras was computed by several
authors~\cite{BLM00,ErdmannHolm99,Locateli99,Zhang97}.

\bibliographystyle{amsplain}
\bibliography{hhcta}

\end{document}